\documentclass{amsart}
\usepackage{amssymb}
\newtheorem{Theo}{Theorem}
\newtheorem{Lem}[Theo]{Lemma}
\newtheorem{Cor}[Theo]{Corollary}
\newtheorem{Prop}[Theo]{Proposition}

\newcommand{\C}{\mathbb{C}}
\begin{document}
\title{Primes in short arithmetic progressions}
\author[J.-C. Schlage-Puchta]{Jan-Christoph Schlage-Puchta}
\maketitle

The Large Sieve inequality in the form 
\[
\sum_{q\leq Q} q\sum_{a=1}^q\left|\sum_{{n\leq N}\atop{n\equiv a\pmod{q}}} a_n
 - \frac{1}{q}\sum_{n\leq N} a_n\right|^2 <(N+Q^2)\sum_{n\leq N} |a_n|^2
\]
is essentially optimal. However, in several applications many of the $a_n$
vanish, and one might expect better estimates then. In fact, such estimates
were given by P. D. T. A. Elliott\cite{Ell}. He showed the following estimate:
\begin{Theo}
$N$ and $Q$ be integers, $a_p$ be complex numbers for all primes $p\leq N$. Then
we have the estimate
\[
\sum_{q\leq Q}(q-1)\sum_{(a,q)=1}\left| \sum_{{p\leq N}\atop{p\equiv
a\pmod{q}}} a_p - \frac{q}{\varphi(q)}\sum_{p\leq N}a_p\right|^2
\ll_\epsilon \left(\frac{N}{\log N} + Q^{54/11+\epsilon}\right)
\sum_{p\leq N}|a_p|^2 
\] 
\end{Theo}
Under GRH, $Q^{54/11}$ may be replaced by $Q^{4+\epsilon}$. In analogy to the
large sieve, he conjectured that one may replace this term by $Q^{2+\epsilon}$.

Using a completely different approach, Y. Motohashi\cite{Mot} showed that
\begin{equation}
\sum\limits_{q\leq Q}{\sum\limits_{\chi\!\!\pmod{q}}}\!\!^*|\pi(x, \chi)|^2\leq
\frac{(2+o(1))x^2}{\log x\log x/Q^{1/2}}
\end{equation}
for $x>Q^{5+\epsilon}$, where
$\pi(x, \chi)=\sum_{p\leq x}\chi(p)$. He also conjectured, that
$Q^{5+\epsilon}$ may be replaced by $Q^{2+\epsilon}$.

Here we will combine the Large Sieve principle with Selberg's sieve to prove
the conjecture of Elliott and give a version of (1) valid for
$x>Q^{2+\epsilon}$. 

I would like to thank D. R. Heath-Brown for his help on Proposition
\ref{CharHilf} which allowed me to reduce the exponent to
$2+\epsilon$, and the referee for pointing out some mistakes.

\begin{Theo} \label{Haupt}
Let $N$ and $Q$ be integers with $N>Q^{2+\epsilon}$, $a_p$ be complex
numbers for any prime $p\leq N$, and let $2\leq R\leq \sqrt{N}$ be an
integer. Then we have the estimate 
\[
\sum_{q\leq Q}{\sum_{\chi\pmod{q}}}^*\left|\sum_{p\leq N}
a_p\chi(p)\right|^2 \ll_\epsilon\frac{N}{\log N}\sum_{p\leq N}|a_p|^2
\]
\end{Theo}

As this estimate is the analogue of the large sieve estimate, we can
give analogues of Hal\'asz-type inequalities, too. As there is a
variety of different large value estimates, the same is true for these
bounds. However, since the optimal estimate depends on the particular
application, we only mention the following: 

\begin{Theo} \label{Halasz}
Let $q$ be an integer. Let $\mathcal{C}$ be a set of characters
$\pmod{q}$, $a_p$ be complex numbers for any prime $p\leq N$. Then we
have for $k=2, 3$ or, if $q$ is cubefree, for any integer $k\geq 2$,
the estimates
\[
\sum_{\chi\in\mathcal{C}}\left|\sum_{p\leq N}a_p\chi(p)\right|^2 \leq
\left(\frac{N}{\log R} + c_{k, \epsilon}N^{1-1/k}q^{(k+1)/(4k^2)+\epsilon}
|\mathcal{C}|R^{2/k}\right)\sum_{p\leq N} |a_p|^2
\]
and
\[
\sum_{\chi\in\mathcal{C}}\left|\sum_{p\leq N}a_p\chi(p)\right|^2 \leq
\left(\frac{N}{\log R} + R^2|\mathcal{C}|\sqrt{q}\log q\right) \sum_{p\leq
N} |a_p|^2. 
\]
If $\mathcal{C}$ is a set of characters to moduli $q\leq Q$, the same
bounds apply with $q$ replaced by $Q^2$, where $k$ can be chosen
arbitrarily, if all occuring values of $q$ are cubefree, and $k=2, 3$
otherwise. 
\end{Theo}

From this we conclude immediately

\begin{Cor} \label{Mot}
We have for $x>Q^{2+\epsilon}$ the estimate
\[
\sum_{q\leq Q}{\sum_{\chi\!\!\pmod{q}}}\!\!^*|\pi(x, \chi)|^2\leq
C_\epsilon\frac{x^2}{\log^2 x}
\]
Moreover, for $x>Q^{3+\epsilon}$ this can be made completely
explicit:
\[
\sum_{q\leq Q}{\sum_{\chi\!\!\pmod{q}}}\!\!^*|\pi(x, \chi)|^2\leq
\frac{(2+o(1))x^2}{\log x\log x/Q^3}
\]
\end{Cor}

We can also consider a single character:

\begin{Cor} \label{einer}
Let $\chi$ be a complex character. Then we have
\[
|\pi(x,\chi)|\leq\left(\left(\frac{1+\phi/\alpha}{2-2\phi/\alpha}\right)^{1/2}
+ o(1)\right)
\frac{x}{\log x},
\]
where $\alpha=\frac{\log x}{\log q}$
and $\phi=\frac{1}{4}$ if $q$ is cubefree, and $\phi=\frac{1}{3}$
otherwise. 
\end{Cor}

Note that this estimate is nontrivial as soon as $x>q^{3/4}$
resp. $x>q$, depending on whether $q$ is cubefree or not.
With a little more work, we obtain the following statement.

\begin{Cor}\label{Rest}
Let $D, x, Q$ be parameters with $x>Q^{1+\epsilon}D^2$. Let $N$ be the number 
of moduli $q\leq Q$, such that there is some primitive character $\chi$ of
order $d\leq D$ and some $d$-th root of unity $\zeta$, such that there is no
prime $p\leq x$ with $\chi(p)=\zeta$. Then we have $N\ll_\epsilon D$.
\end{Cor}

This was proven by Elliott with $D=3$ under the condition
$x>Q^{54/11+\epsilon}$.

We begin the proof of our Theorems with the following two sieve
principles. 

\begin{Lem}[Bombieri]\label{Bom}
Let $V, (\cdot, \cdot)$ be an inner product space, $v_i\in V$. Then for any
$\Phi\in V$ we have
\[
\sum_i |(\Phi, v_i)|^2 \leq \|\Phi\|^2\max_i\sum_j |(v_i, v_j)|
\]
\end{Lem}

This is Lemma 1.5 in \cite{Mon}.

\begin{Lem}[Selberg]\label{Selberg}
Let $R, N$ be integers, such that $R^2<N$. Then there is a function $g$, which
has the following properties:
\begin{enumerate}
\item $g(1)=1$, $|g(n)|\leq 1$ for $n\leq R$, $g(n)=0$ for $n>R$.
\item $\sum_{n\leq N} \big((1\ast g)(n)\big)^2 \leq \frac{N}{\log R} + R^2$
\end{enumerate}
\end{Lem}

This is the usual formulation of Selberg's sieve when used to count
the set of primes $\leq N$, see e.g. \cite{HR}, chapter 3, especially
Theorem 3.3. In the sequel, we will denote the function given by Lemma
\ref{Selberg} with $g$ and set $f=\big(1\ast g\big)^2$. We will have
to bound character sums involving $f$, these computations are
summarized in the following Proposition. 

\begin{Prop}\label{CharHilf}
Let $\chi\pmod{q}$ be a character, $R, N, f$ and $g$ as in Lemma
\ref{Selberg}, and define $S=\sum_{n\leq N} f(n)\chi(n)$.
\begin{enumerate}
\item If $\chi$ is principal, we have $|S|<\frac{N}{\log R}+R^2$.
\item Assume that $\chi$ is nonprincipal. Then we have for any fixed $A$ the
estimate $\sum_{\nu=1}^\infty f(\nu)\chi(\nu) e^{-\log^2(\nu/N)}
\ll_{\epsilon, A} R^2q^{1/2}\left(\frac{N}{R^2q}\right)^{-A}$. 
\item If $\chi$ is nonprincipal, we have the bounds $|S|\leq
R^2\sqrt{q}\log q$ and $|S|\leq c_{k, \epsilon}R^{2/k} N^{1-1/k}
q^{(k+1)/(4k^2)+\epsilon}$ for $k=2, 3$, or, if $q$ is cubefree, for
$k\geq 2$ arbitrary. 
\end{enumerate}
\end{Prop}
{\em Proof:}
The first assertion is already contained in Lemma \ref{Selberg}.

Assume now that $\chi$ is nonprincipal. Then we have
\begin{eqnarray*}
\left|\sum_{n\leq N} f(n)\chi(n)\right| & = & \left|\sum_{n\leq N}
\left(\sum_{d|n} g(d)\right)^2\chi(n)\right|\\
 & = & \left|\sum_{d_1, d_2\leq R} g(d_1)g(d_2)\chi\left([d_1,d_2]\right)
\sum_{n\leq N/[d_1, d_2]} \chi(n)\right|\\ 
 & \leq & \sum_{d_1, d_2\leq N} |g(d_1)g(d_2)| \cdot\left|
\sum_{n\leq N/[d_1, d_2]} \chi(n)\right|\\
 & \leq & \sum_{d_1, d_2\leq R}\left|\sum_{n\leq N/[d_1, d_2]} \chi(n)\right|
\end{eqnarray*}
The inner sum can be estimated using either the Polya-Vinogradoff-inequality
or Burgess estimates, leading to $|S|\leq R^2\sqrt{q}\log q$ resp.
$|S|\leq c_{k, \epsilon} R^{2/k}N^{1-1/k}q^{(k+1)/(4k^2)+\epsilon}$, thus we obtain
the third statement.

To prove the second statement, we begin as above to obtain the
inequality
\[
\left|\sum_{n=1}^\infty f(n)\chi(n)e^{-\log^2(n/N)}\right|\leq
\sum_{d_1, d_2\leq R} \left|\sum_{n=1}^\infty \chi(n)
e^{-\log^2([d_1, d_2]n/N)}\right|
\]
Write $d=[d_1, d_2]$. Using the Mellin-transform
$\frac{1}{2\sqrt{\pi} i}\int_{(2)} x^{-s} e^{s^2/4} ds = e^{-\log^2 x}$,
the inner sum can be expressed as 
\[
\sum_{n=1}^\infty \chi(n)e^{-\log^2(dn/N)} = \frac{1}{2\sqrt{\pi}
i}\int_{(2)}L(s, \chi) e^{s^2/4} (N/d)^s ds
\]
Now we shift the path of integration to the line $\Re s=-A$ with
$A>0$. Denote with $\chi_1$ the primitive character inducing
$\chi$. Then we have
\[
L(s, \chi)=\prod_{p|q_2}\left(1-\chi_1(p)p^{-s}\right)
L(s, \chi_1).
\]
For $A>2$, the first factor is $\ll q_2^A$, whereas the $L$-series can
be estimated using the functional equation to be $\ll
(q_1(|t|+2))^{A+1/2}$, hence the right hand side is $\ll_A 
q^{1/2}\left(\frac{N}{dq}\right)^{-A} \leq q^{1/2}
\left(\frac{N}{R^2 q}\right)^{-A}$. Hence the whole sum can be bounded
by $c(A)R^2q^{1/2}\left(\frac{N}{R^2 q}\right)^{-A}$.  

To prove Theorem \ref{Haupt}, we follow the lines of the proof of the large
sieve resp. the Hal\'asz-inequality, however, we apply Lemma \ref{Bom} to
a different euclidean space. Consider the subspace $V<l^\infty$ consisting of
all bounded sequences $(a_n)$, such that $a_n=0$ whenever $f(n)=0$,
where $f$ is 
defined as in Lemma \ref{Selberg}. On this space define a product as
$\langle (a_n), (b_n)\rangle := \sum_{n=1}^\infty f(n)e^{-\log^2(n/N)}
a_n \overline{b_n}$. 
Now we apply Lemma \ref{Bom} to this space and the set of vectors
$\Phi=(\hat{a}_n)$, where $\hat{a}_p=a_pe^{\log^2 p/N}$, for prime
numbers $p$ in the range $R^2<p\leq N$, and
$\hat{a}_n=0$ otherwise, and $v_i=(\hat{\chi}(n))$, where similary
$\hat{\chi}(n)=\overline{\chi(n)}$, if $f(n)\neq 0$, and 0
otherwise. Now the inequality reads as
\begin{eqnarray*}
\sum_{q\leq Q}\sum_{\chi\pmod{q}}\left|\sum_{R^2<p\leq N}
a_p\chi(p)\right|^2 & \leq & 
\max_{\chi}\left(\sum_{n=1}^\infty f(n)e^{-\log^2 n/N} +
\sum_{\chi'\neq\chi} \left|\sum_{n\leq N}f(n)e^{-\log^2
n/N}\chi\overline{\chi'}(n) \right| \right)\\
&&\quad\times\sum_{p\leq N} |a_p|^2e^{2\log^2(p/N)}
\end{eqnarray*}
where the maximum is taken over all characters with moduli at most $Q$.
From Lemma \ref{Selberg} it follows that the first term inside the
brackets is $\ll\frac{N}{\log R}$, provided that $R<N^{1/3}$, say. For
the second term, let $\chi$ be 
a character $\pmod{q}$ and  $\chi'$ a character
$\pmod{q'}$. Then $\chi\overline{\chi'}$ is a character $\pmod{[q,
q']}$. By Proposition \ref{CharHilf}, each term in the outer sum can
be bounded by $c(A)R^2[q, q']^{1/2}\left(\frac{N}{R^2[q, q']}
\right)^{-A}$, hence the whole sum is $\leq c(A)Q^3R^2
\left(\frac{N}{R^2Q^2}\right)^{-A}$. Since by asumption
$N>Q^{2+\epsilon}$, we can choose $R=Q^{\epsilon/4}$, $A=6/\epsilon+1$
to bound this by some constant depending only on $\epsilon$. Thus we
get the estimate
\[
\sum_{q\leq Q}\sum_{\chi\pmod{q}}\left|\sum_{R^2<p\leq N}
a_p\chi(p)\right|^2 \ll \left(\frac{N}{\epsilon\log N} +
C_\epsilon\right)\sum_{p\leq N} |a_p|^2.
\]
The range $n\leq R^2$ can be estimated using the usual large sieve
inequality, which gives $(R^2+Q^2)\sum_{p\leq N} |a_p|^2$, which is
negligible. Hence Theorem \ref{Haupt} is proven.

The proof of Theorem \ref{Halasz} is similar, but simpler. First,
assume that all characters in $\mathcal{C}$ are characters to a single
modulus $q$. We consider
the vector space $V<\C^N$ consisting of sequences $(a_n)_{n=1}^N$ with 
$a_n=0$ for all $n$ with $f(n)=0$ and the scalar product $\langle
(a_n), (b_n)\rangle := \sum_{n\leq N} f(n) a_n
\overline{b_n}$. Applying Lemma \ref{Bom} as above, we obtain the
estimate
\[
\sum_{\chi\in\mathcal{C}}\left|\sum_{p\leq N}
a_p\chi(p)\right|^2 \leq \left(\frac{N}{\log R} + R^2 +
\big(|\mathcal{C}|-1\big)\Delta(R, N, q)\right) \sum_{R\leq p\leq N}
|a_p|^2 
\]
where $\Delta(R, N, q)$ is the bound obtained by Proposition \ref{CharHilf},
i.e. $\Delta(R, N, q)\leq R^2\sqrt{q}{\log q}$, resp. $\Delta(R, N, q)<
c_{k, \epsilon} q^{(k+1)/(4k^2)+\epsilon}N^{1-1/k}R^{2/k}$. The term
$R^2$ can be neglected in comparison with $\Delta(R, N, q)$. This is
obvious in the first case. In the second case, we may assume that
$\Delta(R, N, q)<N$, since otherwise Theorem \ref{Halasz} is an immediate
consequence of the Cauchy-Schwarz-inequality. This implies
$R<N^{1/2}q^{-(k+1)/(2k)}$, which in turn 
implies $R^2<N^{1-1/k}q^{-1-1/k}<\Delta(R, N, q)$. Hence we obtain
Theorem 3 for sets of characters belonging to a single modulus.

The proof for the case that the characters belong to different moduli
is similar, note that $[q_1, q_2]$ is cubefree, if both $q_1$ and
$q_2$ are cubefree.

In the range $Q^{2+\epsilon}\leq x<Q^{3+\epsilon}$, Corollary
\ref{Mot} follows from Theorem \ref{Haupt} by choosing $a_p=1$ for all
prime numbers $p\leq N$, whereas in the range $x>Q^{3+\epsilon}$ it
follows from Theorem \ref{Halasz}. Similarly we obtain corollary
\ref{einer} from Theorem \ref{Halasz}. We choose $\mathcal{C} =
\{\chi_0, \chi, \overline{\chi}\}$ to obtain the estimate
\[
|\pi(x)|^2 + 2|\pi(x, \chi)|^2\leq \frac{x}{\log
 c_{k, \epsilon} x^{1/2}q^{(k+1)/(8k)+\epsilon}}\pi(x)
\]
and choosing either $k=3$ or $k\nearrow\infty$ we obtain the result by
solving for $|\pi(x, \chi)|$.

To prove corollary \ref{Rest}, let $\mathcal{P}$ be the set of prime numbers $p$,
such that there is some character $\chi$ of order $d$ as described in the
corollary. For every such $p$, choose such a character $\chi_1$ together with
all its
powers, and denote the set of all these character with $\mathcal{C}$. Let $\zeta$
be a $d$-th root of unity. We have
\[
\sum_{\chi^d=\chi_0\atop \chi\neq\chi_0}|\pi(x, \chi)|^2 = d\sum_{a=1}^d
\left|\#\{p\leq x|\chi_1(p)=\zeta^a\} - \frac{1}{d}\pi(x, \chi_0)\right|^2
\]
Since by assumption, one of the terms on the right hand side is large,
the right hand side is $\gg\frac{x^2}{d\log^2 x}\geq \frac{x}{D\log^2
x}$. Now we have $|\mathcal{C}|\leq D\cdot|\mathcal{P}|$, thus we get
\[
|\mathcal{P}| \frac{x^2}{D\log^2 x} \ll \frac{x^2}{\log x\log R} +
xDR^2|\mathcal{P}|Q\log Q 
\]
If $D^2Q\log Q<x^{1-\epsilon}$, we can choose $R=x^{\epsilon/4}$, and
the second term on the right hand side is still of lesser order then
the left hand side. With this choice the inequality can be simplified
to $|\mathcal{P}|\ll_\epsilon D$.

\end{document}